\documentclass[11pt,reprint]{revtex4}

\usepackage{graphicx,epsfig,amssymb,amsbsy,amsmath,amsthm}
\newcommand{\bb}[1]{\left({#1}\right)}					
\newcommand{\sq}[1]{\left[#1\right]}						
\newcommand{\cc}[1]{\left\{#1\right\}}					
\newcommand{\op}[1]{\mathcal{#1}}
\newcommand{\ord}[1]{{\sf O}\bb{#1}}					
\newcommand{\abs}[1]{\left|#1\right|}					

\newcommand{\sfrac}[2]{\mbox{$\frac{#1}{#2}$}}	

\newcommand{\hf}{\mbox{$\frac12$}}

\newcommand{\e}{\tilde\eps}
\newcommand{\f}{{\bf f}}
\newcommand{\F}{{\bf F}}

\newcommand{\sh}{h}

\renewcommand{\k}{{k}}

\newtheorem{example}{Example}

\renewcommand{\Re}{{\sf Re\;}}
\renewcommand{\Im}{{\sf Im\;}}
\newcommand{\sign}{\operatorname{sign}}
\newcommand{\erf}{\operatorname{Erf}}

\newcommand{\arctanh}{\operatorname{arctanh}}
\def\eps{\varepsilon}

\newcommand{\fref}[1]{figure~\ref{#1}}
\newcommand{\Fref}[1]{Figure~\ref{#1}}

\newcommand{\eref}[1]{(\ref{#1})}
\newcommand{\erefs}[2]{(\ref{#1})-(\ref{#2})}
\newcommand{\sref}[1]{section~\ref{#1}}
\newcommand{\srefs}[2]{sections~\ref{#1}-\ref{#2}}

\begin{document}

\title{The ghosts of departed quantities in switches and transitions} 
\author{Mike R. Jeffrey}
\affiliation{Engineering Mathematics, University of Bristol, Merchant Venturer's Building, Bristol BS8 1UB, UK, email: mike.jeffrey@bristol.ac.uk, research supported by EPSRC Fellowship grant EP/J001317/2.}

\date{\today}

\begin{abstract}
Transitions between steady dynamical regimes in diverse applications are often modelled using discontinuities, but doing so introduces problems of uniqueness. No matter how quickly a transition occurs, its inner workings can affect the dynamics of the system significantly. Here we discuss the way transitions can be reduced to discontinuities without trivializing them, by preserving so-called {\it hidden} terms. We review the fundamental methodology, its motivations, and where their study seems to be heading. We derive a prototype for piecewise-smooth models from the asymptotics of systems with rapid transitions, sharpening Filippov's convex combinations by encoding the tails of asymptotic series into nonlinear dependence on a switching parameter. We present a few examples that illustrate the impact of these on our standard picture of smooth or only piecewise-smooth dynamics. 
\end{abstract}

\maketitle

\section{{\it Natura non facit saltus}, so the aphorism goes...}\label{sec:intro}

Whether or not Nature makes jumps, mathematical {\it models} can do. By making jumps, those models may become not only simpler for certain systems, but also a better reflection of our true state of knowledge. Yet fundamental questions remain about the uniqueness of flows with discontinuous vector fields, and whether their {\it non-}uniqueness actually offers physical insight into discontinuities as models of physical behaviour. 
Rigorous ideas from the theories of {piecewise-smooth dynamics} and {singular perturbations} are beginning to shed light on the problem. Here we introduce piecewise-smooth dynamics a little differently to usual and, through some simple examples, show the roles and uses of the ambiguity that accompanies a discontinuity.

Many dynamic systems involve intervals of smooth steady change punctuated by sharp transitions. Some we take for granted, such as light refraction or reflection, electronic switches, and physical collisions. 
In elementary mechanics, for example, when two objects collide, a switch is made between `before' and `after' collision regimes, which are each themselves well understood. The patching of the two regimes leaves certain artefacts, such as the choice between a physical rebound solution, and an unphysical (or virtual) solution in which the objects pass through each other without deviating. 
More exotic switches arise in climate models, for instance as a jump in the Earth's surface albedo at the edge of an ice shelf \cite{hartmann94,engler13}, in superconductivity as a jump in conductivity at the critical temperature \cite{bs08}, in models of cellular mitosis \cite{Tyson}, in dynamics of socio-economic and ecological decision implementation \cite{piltz13,radi13,canadabank}, and so on. 

In the case of the collision model, we do not find the discontinuity or virtual solutions too disturbing when first encountered, and move on to apply such insights to the dynamics of nonlinear mechanical systems, and thereafter to electronics, the climate, living processes, etc., perhaps becoming too comfortable with patching over the joins in our models. 
The models seem to work, but calculus requires continuity, so it seems futile to look deeper. Fortunately, the mathematics of matching such `piecewise-defined' systems turn out to be richer than might be expected. 

Consider a system whose behaviour is modelled by a system of ordinary differential equations ${\bf\dot x}=\f({\bf x};y)$, where $y\approx\sign(\sh({\bf x}))$ for some smooth function $\sh$. Our first aim here is to show that, for many classes of behaviour, such approximations take the form
\begin{equation*}
{\bf\dot x}=\f\bb{{\bf x};y(\sh)}\qquad\mbox{where}\qquad
y(\sh)=\sign(\sh)+\ord{\eps/\sh}\;,
\end{equation*}
for arbitrarily small $\eps$. Our second aim is to show why the tails of these expansions matter, and how they can be retained in a piecewise-smooth model as $\eps\rightarrow0$. 
This information seems not to be part of established piecewise dynamical systems theory, but their omission is easily remedied.

The modern era of piecewise-smooth systems begins with Filippov and contemporaries, who showed that differential equations with ``{\it discontinuous righthand sides}'' can at least be solved (e.g. in \cite{aizerman12,f64,f88}).  
What those solutions look like remains an active and flourishing field of enquiry.

\begin{figure}[h!]\centering\includegraphics[width=0.9\textwidth]{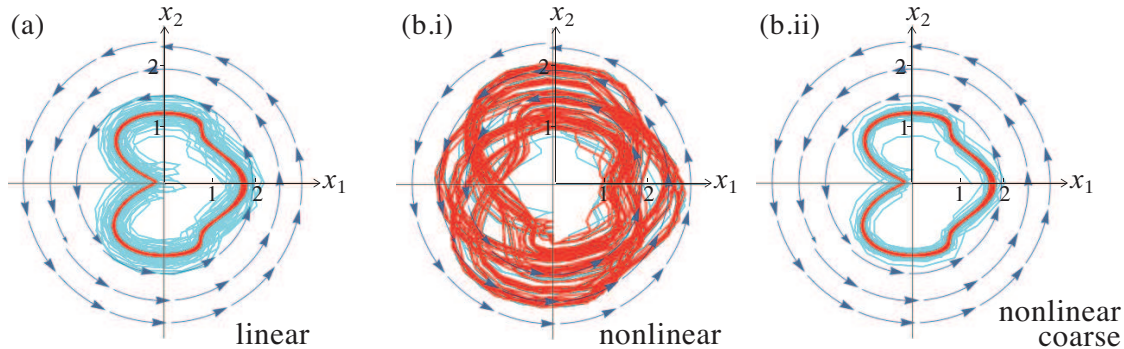}
\vspace{-0.3cm}\caption{\small\sf Dynamics of \eref{whirl}, showing the flow at time $t=0$, and a solution from an initial point $(1,0)$ simulated up to a time $t=2000$ ($t<1000$ shown lightly as transients), by smoothing out the discontinuity with $\lambda=\tanh(x_1/\eps)$, $\eps=10^{-3}$, and 
using explicit Euler discretization in time steps of size $s$, where: (a) $s=10^{-5}$ simulating equation \eref{whirllin} ($s=10^{-4}$ gives similar); (b.i) $s=10^{-5}$ simulating equation \eref{whirlnon}; (b.ii) $s=10^{-4}$ simulating equation \eref{whirlnon}.}\label{fig:whirls}\end{figure}

As an example, take an oscillator given by $\dot x_2=x_1$ and $\dot x_1=-0.01x_1-x_2-\sin(\omega t)$, where the forcing $\sin(\omega t)$ overcomes the damping $-0.01x_1$ to produce sustained oscillations. Say the frequency $\omega$ switches between two values, $\omega=\pi/2$ for $x_1<0$ and $\omega=3\pi/2$ for $x_1>0$. The method usually used to study such switching is due to Filippov \cite{f88,u92,krg03,bc08}, and handles the discontinuity at $x_1=0$ by taking the convex combination of the two alternatives for $\dot x_1$, 
\begin{subequations}\label{whirl}
\begin{equation}\label{whirllin}
\dot x_2=x_1\;,\qquad\dot x_1=-0.01x_1-x_2\;-\;\bb{\sfrac{1+\lambda}2\sin\sq{\sfrac32\pi t}+\sfrac{1-\lambda}2\sin\sq{\hf\pi t}}\;,
\end{equation}
where $\lambda=\sign(x_1)$ for $x_1\neq0$ and $\lambda\in[-1,+1]$ on $x_1=0$. 
We could instead take a convex combination of the frequencies themselves, $\omega=(1+\hf\lambda)\pi$ with $\lambda$ as above, writing
\begin{equation}\label{whirlnon}
\dot x_2=x_1\;,\qquad\dot x_1=-0.01x_1-x_2\;-\;\sin\sq{(1+\hf\lambda)\pi t}\;.\qquad\qquad\quad
\end{equation}
\end{subequations}
\Fref{fig:whirls} shows that the two models have significantly different behaviour. While the linear switching model (a) has a simple limit cycle, the nonlinear model (b.i) has a complex (perhaps chaotic) oscillation. 
This system has been chosen to be deliberately challenging on two counts. 

Firstly, the simulation method matters, particularly to obtain \fref{fig:whirls}(b.i). There are currently no standard numerical simulation codes that can handle discontinuous systems with complete reliability, because although event detection will locate a discontinuity, it is insufficient to determine what dynamics should be applied there. Throughout this paper we show why this question of `what dynamics' should be applied is non-trivial. For reproducibility, \fref{fig:whirls} smooths the discontinuity (replacing the step with a sigmoid function), then uses the Euler method with fixed time step. To find (b.i) requires a numerical simulation with sufficiently precise discretization (see caption), and decreasing precision instead gives (b.ii) (while having no qualitative effect on (a)). 

Secondly, it seems instinctively inconceivable that the systems of equations \eref{whirllin} and \eref{whirlnon} can have different behaviour, because they are identical for $x_1\neq0$, and trajectories cross the zero measure set $x_1=0$ transversally. Despite this, solutions can linger on or even travel along $x_1=0$, whereupon they are influenced differently by linear or nonlinear dependence on $\lambda$, in \eref{whirllin} or \eref{whirlnon} respectively. This behaviour will provide the explanation for \fref{fig:whirls}, to be shown in \sref{sec:whirl}. 

We may then ask whether the behaviour in \fref{fig:whirls}(b.i) is an aberration of the simulation method, or the true behaviour of \eref{whirlnon} as a discontinuous system. 
We may also ask how an observer would interpret this discrepancy if setting up an experiment modelled by \eref{whirlnon}, observing \fref{fig:whirls}(b.i) in experiment, while simulations give \fref{fig:whirls}(a) or (b.ii).  
What we will show is that nonlinear dependence on $\lambda$ introduces fine structure to the switching process, which is captured in (b.i), but can be missed in a less precise simulation as in (b.ii) or by neglecting nonlinearity outright as in (a). 

The consequences of overlooking such nonlinearities of switching can be far more severe. A few key examples are given in \sref{sec:pws}, many more are now appearing in the literature (see e.g. \cite{j15hbifs,j15exit}), but our main aim is to see how they arise and learn how to analyse them.

The starting point to a more general approach to piecewise-smooth systems is quite simple. If a quantity $\f $ switches between states $\f ^+$ and $\f ^-$ as a threshold $\Sigma$ is crossed, $\f $ can be expressed as
\begin{equation}\label{intro}
\f =\hf\bb{1+\lambda}\f ^++\hf\bb{1-\lambda}\f ^-+(\lambda^2-1){\bf g}(\lambda)\;,
\end{equation}
where a step function $\lambda$ switches between $\pm1$ across $\Sigma$ and in $[-1,+1]$ on $\Sigma$. The first two terms have an obvious interpretation, namely the linear interpolation across the jump in $\f $. The last term is less obvious, but the role and origins of each term in \eref{intro} are what we seek to understand here. 

Strictly speaking \eref{intro} may be treated as a differential inclusion. 
When $\bf x$ lies on the switching surface $\Sigma$, the value of $\lambda$ in the set $[-1,+1]$ can usually be fixed by admitting only values of $\f $ that offer viable trajectories: either crossing $\Sigma$ or `sliding' along it. This admissible $\lambda$ value is unique in many situations of interest (but not always at certain singularities, see \cite{j15higen}). The term $\bf g$ is `hidden' outside $\Sigma$ where $\lambda=\pm1$, because its multiplier $\lambda^2-1$ vanishes, but $\bf g$ is potentially crucial inside $\Sigma$ where $\lambda\in[1,+1]$. We shall show how it represents the ``{\it ghost of departed quantities}''~\footnote{a term from {\it The Analyst} \cite{berk1734}} of the transition, called {\it hidden dynamics} \cite{j13error,j15douglas,hairer15}, with significant consequences for local and global behaviour. 

To handle the discontinuity unambiguously we `blow up' the switching {\it surface} $\Sigma$ into a switching {\it layer}, which can reveal 
hidden phenomena such as novel attractors and bifurcations \cite{j15hbifs}. These concepts have been introduced recently, with some heuristic \cite{j13error} and some rigorous \cite{j15douglas} justification. Here we provide a more substantive derivation based on asymptotic transition models. 

We begin in \sref{sec:asydis} by deriving \eref{intro} as a uniform model of switching. The argument begins with a general asymptotic expression of a switch, representative of various differential, integral, or sigmoid models that exhibit abrupt transitions. We delay exploring the motivations for this model to \sref{sec:bi}, as it is somewhat discursive, since discontinuities arise in so many contexts and yet in similar form. 

The mathematics itself in these sections is quite standard, but the universal occurrence of the $\sign$ function and its relation to discontinuous approximations is often under-appreciated, particularly in piecewise-smooth dynamical theory and its applications. Our main aim is to redress this, to show the universality of these expansions for a variety of model classes, and develop the foundations of piecewise-smooth dynamical theory beyond Filippov's convex combinations (but still within Filippov's differential inclusions, see \cite{f88} for both). In \sref{sec:layer} we review how to solve such piecewise-smooth systems. A few stark examples hint at the consequences for piecewise-smooth dynamics 
in \sref{sec:pws}, particularly concerning the passibility of a switching surface, and the novel attractors that may arise. 

To put this more briefly: \sref{sec:asydis} shows how and why nonlinear terms accompany discontinuities, \sref{sec:layer} reviews briefly how to study dynamics at discontinuities, then \sref{sec:pws} combines these to show counterintuitive phenomena caused by such nonlinearity. Finally, \sref{sec:bi} explores some general origins of switching to which the preceding analysis applies, and 
continuing avenues of study are suggested in \sref{sec:conc}.

\section{Asymptotic discontinuity}\label{sec:asydis}

Consider a system that is characterized as having different regimes of behaviour, say
\begin{equation}\label{xasy}
\begin{array}{l}
{\bf\dot x}\sim \f ^+({\bf x})\qquad {\rm for}\quad \sh({\bf x})\gg+\eps\;,\\
{\bf\dot x}\sim \f ^-({\bf x})\qquad {\rm for}\quad \sh({\bf x})\ll-\eps\;,
\end{array}\qquad\qquad\qquad
\end{equation}
where $\f ^+$ and $\f ^-$ are independent vector fields (but each is itself smooth in $\bf x$). 
Some kind of abrupt switch occurs across $|\sh({\bf x})|<\eps$ for small $\eps$. The behaviour inside $|\sh({\bf x})|<\eps$ may be of unknown nature, or of such complexity that our state of knowledge is well represented by the approximation
\begin{equation}\label{xasy0}
\left.\begin{array}{l}
{\bf\dot x}= \f ^+({\bf x})\qquad {\rm for}\quad \sh({\bf x})>0\;\\
{\bf\dot x}= \f ^-({\bf x})\qquad {\rm for}\quad \sh({\bf x})<0\;
\end{array}\right\}\quad{\rm as}\quad\eps\rightarrow0\;.
\end{equation}
The question in either \eref{xasy} or \eref{xasy0} is how to model the system at and around $h=0$. 

For motivation we may consider systems whose full definition we {\it do} know, and which exhibit the behaviour \eref{xasy}-\eref{xasy0}, to derive a common framework for dealing with the switch. 
The result of these investigations (which we delay to \sref{sec:bi} since they are somewhat exploratory), is a typical form near $\sh=0$ which we may represent as a prototype expression
\begin{equation}\label{Fasy}
{\bf\dot x}\;=\;\F ({\bf x},\sh)\;:=\;{\bf p}_0({\bf x})+{\bf p}_1({\bf x})\Lambda(\sh/\eps)+q(\sh/\eps)\sum_{n=1}^\infty{\bf r}_n({\bf x})\;\eps^n/\sh^n 
\end{equation}
in terms of smooth functions ${\bf p_0}({\bf x})$, ${\bf p_1}({\bf x})$, $q(\sh/\eps)$, ${\bf r_n}({\bf x})$, and a sigmoid function
\begin{equation}\label{Leps}
\Lambda(\sh/\eps)\in\left\{\begin{array}{lll}\sign(\sh)&\rm if&|\sh|>\eps\\\sq{-1,+1}&\rm if&|\sh|\le\eps\end{array}\right\}+\ord{\eps/\sh}\;,\qquad \Lambda'(\sh/\eps)>0\;,
\end{equation}
which tends to a discontinuous function, $\Lambda(\sh/\eps)\rightarrow\sign(\sh)$, as $\eps\rightarrow0$. The term ${\bf p}_1\Lambda$ in \eref{Fasy} encapsulates the switching in the system, the summation term contains behaviour that is asymptotically vanishing away from the switch, and the term ${\bf p}_0$ is switch independent. 

The expression \eref{Fasy} is the starting point for the analysis which follows, hereon until \sref{sec:pws}. 

We begin by re-writing \eref{Fasy} in a form that behaves uniformly as $\eps\rightarrow0$. 
Since $\Lambda(\sh/\eps)$ is monotonic in $h$, it has an inverse $V$ such that $h=\eps V(\Lambda)$, then
\begin{equation}\label{fV}
{\bf\dot x}\;=\;\f ({\bf x};\Lambda)\;:=\; {\bf p}_0({\bf x})+{\bf p}_1({\bf x})\Lambda+q(V(\Lambda))\sum_{n=1}^\infty {\bf r}_n({\bf x})(V(\Lambda))^{-n}\;.
\end{equation}
Since this is now a function of $\bf x$ and $\Lambda$, assume that the righthand side of \eref{fV} can be expressed as a formal series in $\Lambda$, 
\begin{equation}\label{fsumgen}
\f({\bf x};\Lambda)=
\sum_{n=0}^\infty{\bf c}_n({\bf x})\Lambda^n\;.
\end{equation}
We can relate the ${\bf c}_n$'s to the ${\bf r}_n$'s, but more useful is to relate them directly to the large $\sh/\eps$ behaviour of ${\bf\dot x}$ in \erefs{xasy}{xasy0}, giving $\f({\bf x},\pm1)\equiv\f ^\pm({\bf x})$. Taking the sum and difference of these gives
\begin{align*}
\hf\bb{\f^+({\bf x})+\f^-({\bf x})}&={\bf c}_0({\bf x})+\sum_{n=1}^\infty{\bf c}_{2n}({\bf x})\;,\\
\hf\bb{\f^+({\bf x})-\f^-({\bf x})}&={\bf c}_1({\bf x})+\sum_{n=1}^\infty{\bf c}_{2n+1}({\bf x})\;,
\end{align*}
so we can eliminate the first two coefficients ${\bf c}_0({\bf x})$ and ${\bf c}_1({\bf x})$ in \eref{fsumgen} to give
\begin{equation}\label{Fgen0}
\f({\bf x};\Lambda)=\frac{\f^+({\bf x})+\f^-({\bf x})}2+\frac{\f^+({\bf x})-\f^-({\bf x})}2\Lambda+\Gamma({\bf x};\Lambda)\;,
\end{equation}
with a remainder term  
\begin{equation}\Gamma({\bf x};\Lambda):=\sum_{n=1}^\infty\cc{{\bf c}_{2n}({\bf x})+\Lambda{\bf c}_{2n+1}({\bf x})}\cc{\Lambda^{2n}-1}\;.
\end{equation}
The factor $\Lambda^{2n}-1$ implies $\Gamma({\bf x};\Lambda)=0$ when $\Lambda=\pm1$. These are the `ghosts' of switching, terms that are lost if we consider only the $\Lambda=\pm1$ states, now retained in the expression $\Gamma({\bf x};\Lambda)$. 
We can take out a factor $\Lambda^2-1$, to find $\Gamma({\bf x};\Lambda)=(\Lambda^2-1){\bf g}({\bf x};\Lambda)$ where
\begin{equation}
{\bf g}({\bf x};\Lambda)=\sum_{n=1}^\infty\sum_{m=0}^{n-1}\cc{{\bf c}_{2n}({\bf x})+\Lambda{\bf c}_{2n+1}({\bf x})}\Lambda^{2m}\;.
\end{equation}

The remaining coefficients ${\bf c}_{n\ge2}$ can in principle be determined from any deeper knowledge we have of $\F$, such as the partial derivatives of $\F$ with respect to $h$ for large or small $\sh/\eps$. For example, if we know the value of $\F({\bf x},0)$ in \eref{Fasy}, then ${\bf c}_{2}=\frac{\f^++\f^-}2-\F({\bf x},0)-\sum_{n=2}^\infty{\bf c}_{2n}$, and successive coefficients can be eliminated by partial derivatives of $\F$ with respect to $h$. In cases where this is not possible, we can propose forms of $\bf g$ based on dynamical or physical considerations (much as we do when proposing dynamical models for smooth systems that may be nonlinear in a state $\bf x$).

The result is that, given an asymptotic expression \eref{Fasy} for a switch across an $\eps$-width boundary in a dynamical system, we obtain an $\eps$-independent form
\begin{eqnarray}\label{Fgen}
{\bf\dot x}=\f ({\bf x};\lambda)= \frac{\f^+({\bf x})+\f^-({\bf x})}2+\frac{\f^+({\bf x})-\f^-({\bf x})}2\lambda+\bb{\lambda^2-1}{\bf g}({\bf x};\lambda)\;,
\end{eqnarray}
as promised in \eref{intro}. This remains valid as $\eps\rightarrow0$, and the switch, whether smooth or discontinuous, is hidden implicitly inside $\lambda$. If we let $\eps\rightarrow0$, then by \eref{Leps} the {\it switching multiplier} $\lambda$ obeys 
\begin{equation}\label{lam}
\lambda\in\left\{\begin{array}{lll}\sign(\sh)&\rm if&|\sh|\neq0\\\sq{-1,+1}&\rm if&|\sh|=0\end{array}\right\}\;.
\end{equation}
The essential point is that we are left with the {\it ghosts}, in the term $\bb{\lambda^2-1}{\bf g}({\bf x}$ which vanishes away from the switch where $\lambda=\pm1$, but does not vanish on $h=0$. The next two sections show how to handle them, and why their existence is non-trivial. 

In \sref{sec:bi} we return to how and when such switches arise in various contexts, including sigmoid-like transitions, higher dimensional ordinary or partial differential equations, and oscillatory integrals. We now turn to the methods used to solve the piecewise-smooth system \eref{Fgen}-\eref{lam}.

\section{The switching layer}\label{sec:layer}

We have derived in \eref{Fgen} an expression for the vector field in our system ${\bf\dot x}=\f ({\bf x};\lambda)$ which, with \eref{lam}, remains valid in the discontinuous limit $\eps\rightarrow0$. One last thing is needed to complete the description of the piecewise-smooth system, and that is to deal with the set-valuedness of $\lambda(0)$ in \eref{lam}. To do this we derive the dynamics on $\lambda$ from the asymptotic relations above. We then derive key manifolds organizing the flow, and interpret the result as a singular perturbation problem. (This section is essentially a review of concepts introduced in \cite{j13error,j15higen}, a modern extension of Filippov's theory \cite{f88}). 

Differentiating $\lambda=\Lambda(\sh/\eps)$ with respect to $t$ gives $\dot\lambda=\Lambda'(\sh/\eps){\dot\sh}/\eps$. Define $\e(\lambda,\eps)=\eps/\Lambda'(\sh/\eps)$, then apply the chain rule and \eref{Fgen} to substitute $\dot\sh={\bf\dot x}\cdot\nabla\sh=\f \cdot\nabla\sh$. Thus at $h=0$ the dynamics of $\lambda$ is given by
\begin{equation}\label{ldash0}
\e\dot\lambda=\f({\bf x};\lambda)\cdot\nabla \sh({\bf x})\qquad{\rm as}\;\;\e\rightarrow0\;.
\end{equation} 
This result is derived in greater detail in \cite{j15higen}, showing that since $\Lambda$ is monotonic by \eref{Leps}, so $\e\ge0$ and $\e\rightarrow0$ as $\eps\rightarrow0$. Since only the limit $\eps\rightarrow0$ concerns us in a piecewise-smooth model, the fact that $\e$ is a function rather than a fixed parameter is of no interest, it is just an infinitesimal like $\eps$.

When we now combine the $\lambda$ dynamics \eref{ldash0} with the original system ${\bf\dot x}=\f({\bf x};\lambda)$, the result is a two timescale system on the switching surface $\sh=0$. Choosing coordinates ${\bf x}=(x_1,x_2,...,x_n)$ where $x_1=\sh({\bf x})$, and writing $\f=(f_1,f_2,...,f_n)$, putting \eref{ldash0} together with \eref{Fgen} we have
\begin{equation}\label{blowup}
\left.\begin{array}{rcl}
\e
\dot\lambda&\!\!=\!\!&f_1(0,x_2,x_3,...,x_n;\lambda)\\
(\dot x_2,...,\dot x_n)&\!\!=\!\!&\bb{f_2(0,x_2,...,x_n;\lambda),...,f_n(0,x_2,...,x_n;\lambda)}
\end{array}\right\}\;\;\;{\rm on}\quad x_1=0\;.
\end{equation}
This defines dynamics inside the switching layer $\cc{\lambda\in\sq{-1,+1},\;(x_2,...,x_n)\in\mathbb R^{n-1}}$, and we call \eref{blowup} the {\it switching layer system}. 

Rescaling time in \eref{blowup} to $\tau=t/\e$, then setting $\e=0$, gives the fast critical subsystem (denoting the derivative with respect to $\tau$ by a prime)
\begin{equation}\label{lfast}
\left.\begin{array}{rcl}
\lambda'&\!\!=\!\!&f_1(0,x_2,...,x_n;\lambda)\\
(x_2',...,x_n')&\!\!=\!\!&\bb{0,\;...,0}
\end{array}\right\}\;\;\qquad\qquad\qquad\qquad\qquad{\rm on}\quad x_1=0\;,
\end{equation}
which gives the fast dynamics of transition through the switching layer. 
The equilibria of this one-dimensional system form the so-called {\it sliding manifold}
\begin{equation}\label{M}
\op M=\cc{(\lambda,x_2,...,x_n)\in[-1,+1]\times\mathbb R^{n-1}:\;f_1(0,x_2,...,x_n;\lambda)=0}\;.
\end{equation}
When $\op M$ exists, it forms an invariant manifold of \eref{blowup} in the $\e=0$ limit, at least everywhere that $\op M$ is normally hyperbolic, which excludes the set where $\frac{\partial f_1}{\partial \lambda}=0$, namely
\begin{equation}\label{L}
\op L=\cc{(\lambda,x_2,...,x_n)\in\op M:\;\frac{\partial\;}{\partial \lambda}f_1(0,x_2,...,x_n;\lambda)=0}\;.
\end{equation}
Isolating the slow system in \eref{blowup}, and setting $\e=0$ in \eref{blowup}, gives the slow critical subsystem
\begin{equation}\label{blowups}
\left.\begin{array}{rcl}
0&\!\!=\!\!&f_1(0,x_2,...,x_n;\lambda)\\
(\dot x_2,...,\dot x_n)&\!\!=\!\!&\bb{f_2(0,x_2,...,x_n;\lambda),...,\;f_n(0,x_2,...,x_n;\lambda)}
\end{array}\right\}\;\;\;{\rm on}\quad x_1=0\;,
\end{equation}
which gives the dynamics on $\op M$ in the limit $\e=0$, called a {\it sliding mode}. 

These elements \eref{Fgen} with \eref{ldash0}, implying \eref{blowup}-\eref{M}, form the basis of everything that follows. We shall look at some of the behaviours they give rise to, hinting at the zoo of singularities and nonlinear phenomena that remain a large classification task for future work. 

In the context of piecewise-smooth systems, we are concerned with these results only in the limit $\e\rightarrow0$, not the perturbation to $\e>0$ that is typically of interest in singular perturbation studies. However, for more general interest it is worth relating these to singular perturbation theory. 
The system \eref{blowup} is the singular limit of
\begin{equation}\label{ublowup}
\left.\begin{array}{rcl}
\e(\lambda,\eps)\dot \lambda&\!\!=\!\!&f_1(\eps u,x_2,x_3,...,x_n;\lambda)\\
(\dot x_2,...,\dot x_n)&\!\!=\!\!&\bb{f_2(\eps u,x_2,...,x_n;\lambda),...,f_n(\eps u,x_2,...,x_n;\lambda)}
\end{array}\right\}\;\;\;{\rm for}\quad \eps\ll1\qquad\quad
\end{equation}
where $u=x_1/\eps$ and $\lambda=\Lambda(x_1/\eps)=\Lambda(u)$ with $\eps\ge0$. Equivalently we can write
\begin{equation*}
\left.\begin{array}{rcl}
\eps\dot u&\!\!=\!\!&f_1(\eps u,x_2,x_3,...,x_n;\Lambda(u))\\
(\dot x_2,...,\dot x_n)&\!\!=\!\!&\bb{f_2(\eps u,x_2,...,x_n;\Lambda(u)),...,f_n(\eps u,x_2,...,x_n;\Lambda(u))}
\end{array}\right\}\;\;\;{\rm for}\quad \eps\ll1\;,\;
\end{equation*}
which is a more commonly seen expression in recent singular perturbation studies of piecewise-smooth systems (see e.g. \cite{ts11}). In \cite{j15higen} it is shown than \eref{blowup} has equivalent slow-fast dynamics to \eref{ublowup} on the discontinuity set $x_1=0$ in the critical limit $\eps=0$.

With this we depart the smooth world. In \sref{sec:asydis} we showed how a prototype asymptotic expansion \eref{Fasy} could be represented as a discontinuous system in the small $\eps$ limit, but left behind nonlinearities in the switching multiplier $\lambda$. We now take our expressions \eref{Fgen}-\eref{lam} valid for $\eps\rightarrow0$, and the dynamics of $\lambda$ at the discontinuity given by \eref{ldash0}, and continue henceforth in the realm of piecewise-smooth dynamics alone (where $\eps=\e=0$), to show some of the counterintuitive phenomena that nonlinear switching terms give rise to.

\section{Hidden dynamics: examples}\label{sec:pws}

To summarize the analysis above: we have a general expression for a discontinuous system in the form \eref{Fgen}-\eref{lam}, for some smooth vector functions $\f ^+$, $\f ^-$, and $\bf g$. Only $\f ^\pm$ are fixed by the dynamics in $\sh\neq0$, with $\bf g$ being directly observable only on $h=0$. 
%
On $h=0$ we look inside the switching layer $\lambda\in[-1,+1]$, whose dynamics is given by the two timescale system \eref{blowup} 
in coordinates where $h=x_1$. If $\dot\lambda=0$ then solutions can become trapped inside the layer, on the sliding manifold $\op M$ if/where it exists, upon which sliding dynamics \eref{blowups} occurs. 
In this section we can replace $\e$ by simply $\eps$, and only the limit $\eps\rightarrow0$ concerns us.

\subsection{Cross or Stick?}\label{sec:h}

Consider what happens when the flow of \eref{Fgen} arrives at a switching surface $\sh=0$. At least one of the vector fields $\f^\pm({\bf x})$ (the one the flow arrived through) points towards the surface. Whether or not the flow then crosses the surface is determined first by the vector field on the other side of the surface, $\f ^\mp({\bf x})$, and possibly also by $\bf g$. 

\begin{figure}[h!]\centering\includegraphics[width=0.75\textwidth]{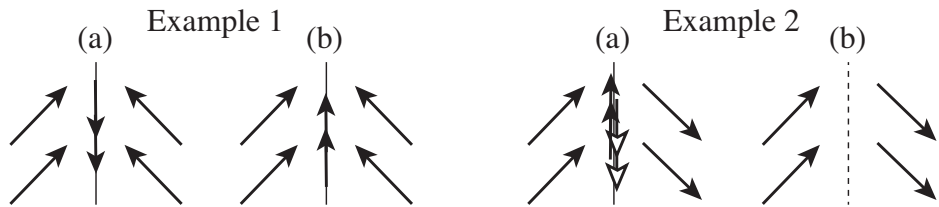}
\vspace{-0.3cm}\caption{\small\sf Sketch of two planar piecewise constant systems. Each portrait (a) includes hidden dynamics that is not obvious outside the switching surface: sliding downwards in Example 1, and in Example 2 an attracting upwards sliding solution and repelling downwards sliding solution. Each portrait (b) excludes hidden dynamics and shows the Filippov dynamics: sliding upwards in Example 1 and crossing in Example 2. }\label{fig:eg1}\end{figure}

If $\f^+({\bf x})\cdot\nabla\sh<0<\f^-({\bf x})\cdot\nabla\sh$ at $\sh=0$, as in \fref{fig:eg1} (Example 1), both vector fields point towards the switching surface so the flow obviously cannot cross it. The normal component $\f({\bf x};\lambda)\cdot\nabla\sh$ changes sign as $\lambda$ changes between $\lambda=\pm1$, so there must exist at least one value $\lambda\in[-1,+1]$ for which $\f({\bf x};\lambda)\cdot\nabla \sh=0$. This defines a so-called {\it sliding mode} on the switching surface, i.e. a solution evolving according to \eref{blowups}. 

If \eref{Fgen} depends linearly on $\lambda$, i.e. if ${\bf g}\equiv0$, then the sliding mode given by \eref{blowups} is unique (and is exactly that described by Filippov \cite{f88}). If ${\bf g}$ is nonzero then there may be multiple sliding modes, and the precise dynamics must be found using \eref{blowup}. 

\begin{example}\label{eg:slide}
A simple example of hidden dynamics is given by comparing the two systems 
$${\rm(a)}\quad(\dot x_1,\dot x_2)=(-\lambda,2\lambda^2-1)\;,\qquad{\rm(b)}\quad(\dot x_1,\dot x_2)=(-\lambda,1)\;,$$
with $\lambda=\sign(x_1)$, shown in \fref{fig:eg1}. These appear to be identical for $x_1\neq0$, where $(\dot x_1,\dot x_2)=(-\sign(x_1),1)$. It is only on $x_1=0$ that their behaviour may differ. To find this we blow up $x_1=0$ into the switching layer $\lambda\in[-1,+1]$, given by applying \eref{blowup}, 
$${\rm(a)}\quad(\eps\dot\lambda,\dot x_2)=(-\lambda,2\lambda^2-1)\;,\qquad{\rm(b)}\quad(\eps\dot\lambda,\dot x_2)=(-\lambda,1)\;,$$
respectively, for $\eps\rightarrow0$. We seek sliding modes by solving $\dot\lambda=0$. Both have sliding manifolds $\op M$ at $\lambda=0$, and therefore sliding modes with, however, contradictory vector fields
$${\rm(a)}\quad(\eps\dot\lambda,\dot x_2)=(0,-1)\;,\qquad\qquad\quad{\rm(b)}\quad(\eps\dot\lambda,\dot x_2)=(0,+1)\;.$$
\end{example}

Hence systems that appear the same outside the switching surface can have distinct, and even directly opposing, sliding dynamics on the surface, due to nonlinear dependence on $\lambda$. 

If $\bb{\f^+({\bf x})\cdot\nabla\sh}\bb{\f^-({\bf x})\cdot\nabla\sh}>0$ at $\sh=0$, as in \fref{fig:eg1} (Example 2), both vector fields point the same way through the switching surface, so the flow may be expected to cross it. If ${\bf g}\equiv0$, in fact, the flow {\it will} cross the surface, because $\f^\pm({\bf x})\cdot\nabla\sh$ have the same sign as each other, so the linear interpolation \eref{Fgen} (as $\lambda$ varies between $\pm1$) cannot pass through zero and there can be no sliding modes (no solutions of \eref{blowups}). If $\bf g$ is nonzero then the flow may stick to the surface, and solutions sliding along the surface are found using \eref{blowups}.

\begin{example}\label{eg:cross}
Taking again $\lambda=\sign(x_1)$, consider the second system in \fref{fig:eg1}, given by
$${\rm(a)}\quad(\dot x_1,\dot x_2)=(2\lambda^2-1,-\lambda)\;,\qquad{\rm(b)}\quad(\dot x_1,\dot x_2)=(1,-\lambda)\;.$$
These both appear the same, $(\dot x_1,\dot x_2)=(1,-\lambda)$, for $x_1\neq0$. 
In the switching layer \eref{blowup} gives 
$${\rm(a)}\quad(\eps\dot\lambda,\dot x_2)=(2\lambda^2-1,-\lambda)\;,\qquad{\rm(b)}\quad(\eps\dot\lambda,\dot x_2)=(1,-\lambda)\;.$$
Solving $\dot\lambda=0$ gives sliding modes $\lambda=\pm1/\sqrt2$ in (a), but no solutions in (b). In (a) the derivative $\eps\partial\dot\lambda/\partial\lambda=\pm2\sqrt2$ reveals that the solutions $\lambda=-1/\sqrt2$ and $\lambda=+1/\sqrt2$ are attracting and repelling respectively, so solutions collapse to $\lambda=-1/\sqrt2$ and follow the sliding dynamics $(\eps\dot\lambda,\dot x_2)=(0,+1/\sqrt2)$. In (b) the fast subsystem $\eps\dot\lambda=1$ carries the solution across the switching layer without stopping. 
\end{example}

Hence even the simple matter of whether or not a system will cross through a switch cannot be determined without considering the effects of nonlinearity at the switch. 

We refer to the behaviour in (a) for each example as `hidden dynamics', because it arises through the addition of hidden terms, $2(\lambda^2-1)$ in both examples, to the linear system (b). 
In Example 2 we could even replace the second component with $\dot x_2=1$ for both (a) and (b), then $(\dot x_1,\dot x_2)=(1,1)$ for $x_1\neq0$, and the discontinuity is an effect localized entirely to $x_1=0$). 

\subsection{Hidden van der Pol system}

\begin{figure}[h!]\centering\includegraphics[width=0.7\textwidth]{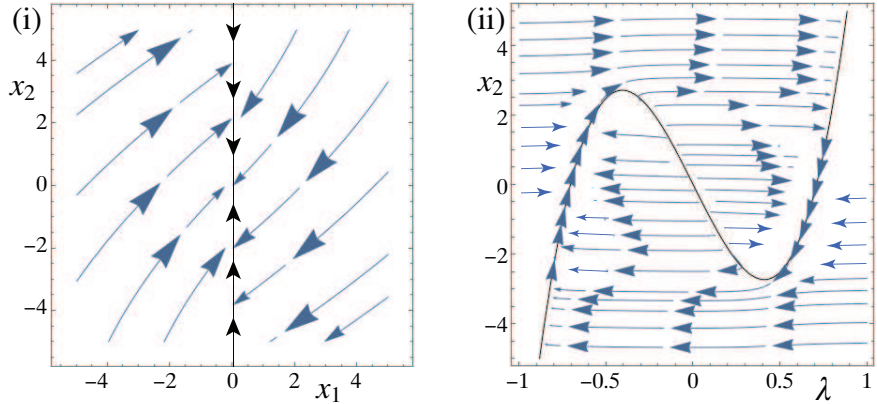}
\vspace{-0.3cm}\caption{\small\sf Simulations of \eref{vdp1} showing: (i) the flow in the $(x_1,x_2)$ plane, (ii) the flow inside $x_1=0$ given by \eref{vdp2}. }\label{fig:vdp}\end{figure}

Hidden dynamics can be much more interesting. 
Take for example the system
\begin{equation}\label{vdp1}
(\dot x_1,\dot x_2)=\bb{\;\sfrac1{10}x_2+\lambda-2\lambda^3,\;-\lambda\;}\qquad\qquad\;\;
\end{equation}
where $\lambda=\sign(x_1)$. This is deceptively simple for $x_1\neq0$, where
\begin{equation}\label{vdp1o}
(\dot x_1,\dot x_2)
=\left\{\begin{array}{lll}\bb{\;\sfrac1{10}x_2-1,\;-1\;}&\rm if&x_1>0\;,\\\bb{\;\sfrac1{10}x_2+1,\;+1\;}&\rm if&x_1<0\;,\end{array}\right.
\end{equation}
illustrated in \fref{fig:vdp}(i). The surface $x_1=0$ is attracting. The switching layer from \eref{blowup}, however, reveals a van der Pol oscillator, 
\begin{equation}\label{vdp2}
(\eps\dot\lambda,\dot x_2)=\bb{\;\sfrac1{10}x_2+\lambda-2\lambda^3,\;-\lambda\;}\;.\qquad\qquad
\end{equation}
To identify the hidden term notice that $\lambda^3=\lambda+(\lambda^2-1)\lambda$, which looks like $\lambda$ for $x_1\neq0$. In Filippov's method we ignore the hidden term $(\lambda^2-1)\lambda$ which vanishes outside $x_1=0$, and then we would find that the point $x_1=x_2=0$ is attracting. Including the nonlinear term, however, the switching variable $\lambda\in[-1,+1]$  undergoes relaxation oscillations hidden inside $x_1=0$. The dynamics inside the switching layer is shown in \fref{fig:vdp}(ii). 

\begin{figure}[h!]\centering\includegraphics[width=0.7\textwidth]{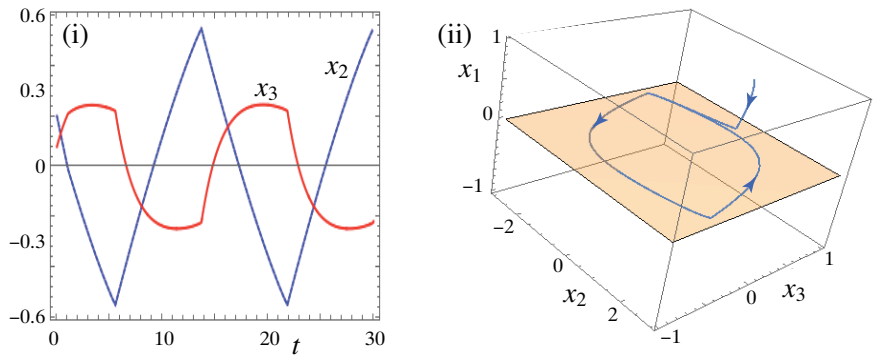}
\vspace{-0.3cm}\caption{\small\sf Simulations revealing the hidden dynamics of \eref{vdp1}: (i) graphs of the variable $x_2$, and of $x_3$ using \eref{vdp3} with $\beta=10^{-4}$, (ii) the corresponding orbit in the space of $(x_1,x_2,x_3)$
, with the switching surface at $x_1=0$. }\label{fig:vdp3d}\end{figure}

The oscillations can be made observable if we plot $x_2$ against time, \fref{fig:vdp3d}(i). Or we can view the dynamics of $\lambda$ itself by coupling the system to a third variable, say
\begin{equation}\label{vdp3}
\beta\dot x_3=\lambda-x_3\;,
\end{equation}
for small $\beta$. A simulation is shown in \fref{fig:vdp3d}(i), with the orbit in phase space shown in (ii). 

Similarly to \fref{fig:whirls}, simulating the hidden dynamics (the oscillation here) is reliant on a sufficiently precise numerical simulation. If we solve by letting $\lambda=\tanh(x_1/\eps)$ with $\eps=10^{-3}$, using an Euler discretization with fixed time step less than or equal to $\eps$ (or using some other more precise method), we obtain \fref{fig:vdp3d}. A more coarse simulation may miss the hidden oscillation, for example with a fixed discretization time step $s\ge4\eps$ the state $x_2$ seems to instead reach the equilibrium $x_1=x_2=0$ of the linear theory (simulations not shown). We will comment more on the general principles behind such sensitivity at the end of \sref{sec:whirl}.

\subsection{Oscillator revisited: hidden dynamics and its robustness}\label{sec:whirl}

Let us extract the hidden term for the oscillator introduced in \eref{whirlnon}. We can write
\begin{eqnarray}
\sin\bb{(1+\hf\lambda)\pi t}&=&\sin\bb{\sfrac{1+\lambda}2\sfrac32\pi t+\sfrac{1-\lambda}2\hf\pi t}\nonumber\\
&=&\sfrac{1+\lambda}2\sin\bb{\sfrac32\pi t}+\sfrac{1-\lambda}2\cos\bb{\hf\pi t}+\bb{\lambda^2-1}g(t;\lambda)
\end{eqnarray}
where some lengthy algebra yields
\begin{eqnarray*}
g(t;\lambda)&=&\sfrac14\bb{\sin[\sfrac32\pi t]+\sin[\hf\pi t]-2\pi t}+\sfrac14
\sum_{i=1}^\infty\bb{\sfrac12\pi t}^{2i+1}\times\\&&\left[\sum_{j=1}^{2i-1}\sfrac1{(2i+1)!}\cc{\bb{\sfrac{1+\lambda}2}^{2j}3^{2i+1}+\bb{\sfrac{1-\lambda}2}^{2j}}-\right.\\&&\left.
\sum_{j=0}^\infty\sfrac{\bb{\pi t/2}^{2j}}{(2i)!(2j+1)!}\cc{\bb{\sfrac{1-\lambda}2}^{2i-1}\bb{\sfrac{1+\lambda}2}^{2j}3^{2j+1}+\bb{\sfrac{1+\lambda}2}^{2i-1}\bb{\sfrac{1-\lambda}2}^{2j}3^{2i}}\right]\;.
\end{eqnarray*}

\begin{figure}[h!]\centering\includegraphics[width=0.9\textwidth]{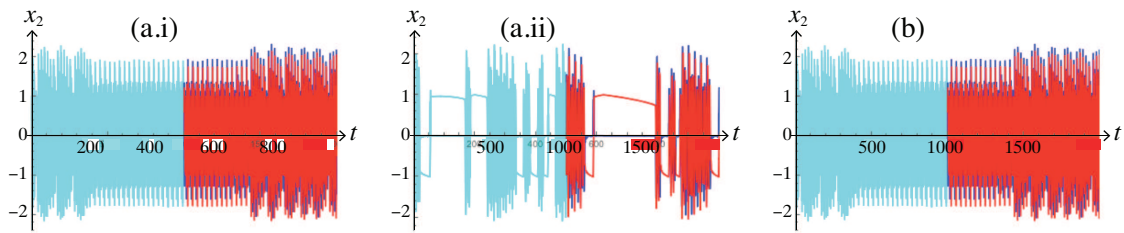}
\vspace{-0.3cm}\caption{\small\sf Simulation of $x_2(t)$ corresponding to \fref{fig:whirls}. Segments of sliding can be seen in (a.i).}\label{fig:whirlt}\end{figure}

The direct effect of the nonlinear term is fairly benign compared to the examples above: it merely slows the dynamics as it crosses the switching surface. The nonlinearity in $\lambda$ means that small regions of sliding, where $\dot\lambda=0$, are able to appear and disappear at $x_1=0$, temporarily preventing solutions from crossing $x_1=0$. They arise from nonlinear terms as in Example 2 above. The sliding can be seen in the simulation of the $x_2=0$ coordinate in \fref{fig:whirlt}. 

\begin{figure}[h!]\centering\includegraphics[width=0.9\textwidth]{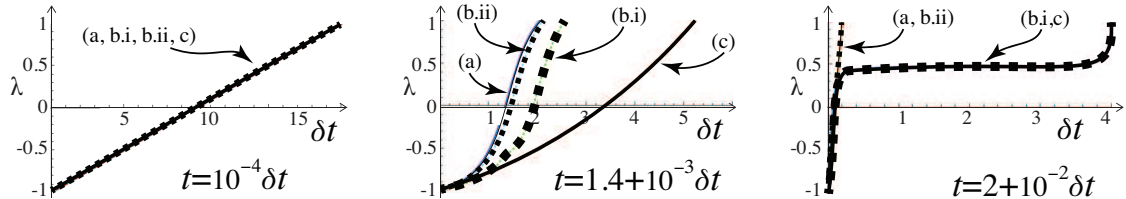}
\vspace{-0.3cm}\caption{\small\sf Simulation of trajectories crossing the layer $\lambda\in[-1,+1]$ in \fref{fig:whirls}, plotting $\lambda$ against the time $\delta t$ spent in transit, at different times $t$ indicated on the figures. Curves are labelled corresponding to \fref{fig:whirls}: (a) linear switching, (b.i) nonlinear switching with fine discretization, (b.ii) nonlinear switching with coarse discretization. Curve (c) shows nonlinear switching with adaptive discretization (using Mathematica's {\sf NDSolve} in default mode). }\label{fig:whirlstime}\end{figure}

In a smoothed-out simulation like \fref{fig:whirls}, this slowing reveals itself as a slowing of trajectories as they attempt to cross $x_1=0$. In \fref{fig:whirlstime} we show this slowing. 
Using the simulation methods described in \fref{fig:whirls}, each graph simulates the evolution of $\lambda$ through the switching layer, and while at $t=0$ all simulations agree, at a later time $t=1.4$ the graph depends strongly on linearity of the model and numerical precision, and at $t=2$ the linear system or coarse simulation are clearly distinct from the nonlinear system, crossing the layer $\lambda\in[-1,+1]$ in much shorter time. This is enough, given the time-dependent sinusoidal control, to alter the connection between trajectories either side of the switch sufficiently and destabilize the oscillation. In the ideal $\eps\rightarrow0$ limit where the switch is discontinuous, this time-lag remains (but $x_1$ remains exactly `sliding' on $x_1=0$ during the switch, rather than slowly transitioning through $|x_1|<\eps$).


It is worth summarizing one result concerning hidden dynamics proposed in \cite{j13error}, where a heuristic case was made that `unmodelled errors' could kill off hidden dynamics, i.e. mask (or essentially eliminate) the nonlinear dependence on $\lambda$ in \eref{Fgen}. Unmodelled errors might include the discretization step of a simulation, time delay or hysteresis of a switching process, or external noise. Essentially, large perturbations by unmodelled errors can kick a system far enough that nonlinear features are missed. We saw in \fref{fig:whirls}(b) how coarse numerical integration killed off the destabilizing effects of nonlinearity. Another simple example would be Example 2 of \sref{sec:h}, where the attracting and repelling sliding modes could be masked in a system with large unmodelled errors (e.g. a discretization step or additive noise larger than the separation between the attracting and repelling branches), so solutions cross as in the linearized $\lambda$ model. 

A similar result in \cite{j13error} showed how, in a dry-friction inspired model, hidden terms can model static friction, but stochastic perturbations of sufficient size destroy it. The outcome was that static and kinetic friction coefficients become equal in more irregular systems. The result was shown rigorously in the presence of white noise in \cite{js2013}. We now summarise the general but partly heuristic result, hoping that the challenge of generalising it will be taken up by future researchers. 

The idea is to add a stochastic perturbation $d{\bf W}$ in the form
\begin{equation}
d{\bf x}=\f({\bf x};\Lambda\bb{\sh(\bf x)/\eps})
dt+\kappa d{\bf W}(t)
\end{equation}
with $\f({\bf x};\lambda)$ given by \eref{Fgen}, with $\Lambda\bb{\sh(\bf x)/\eps}$ being a smooth (or at least continuous) sigmoid function, and $d\bf W$ a standard vector-valued Brownian motion. The zeros of $\f({\bf x};\lambda)\cdot\nabla\sh$ show up as potential wells, maxima and minima of a potential function $U(\sh)=-\int_0^vdu\;\f({\bf x};\Lambda\bb{\sh(\bf x)}/\eps)\cdot\nabla\sh$, which form stationary points of the transition. These correspond to attractors or repellers of the dynamics near $\sh=0$, upon which solutions slide along $\sh=0$. The results of \cite{simpson12} then show that the average state $\langle{\bf x}\rangle$ evolving along $\sh\approx0$ behaves as
 \begin{eqnarray}
\frac{d\langle{\bf x}\rangle}{dt}&=&\f({\bf x};\Lambda\bb{\sh(\bf x)}/\eps)+\ord{\kappa^2}
\end{eqnarray}
recalling \eref{Fgen} 
and $(1-\Lambda^2\bb{\sh(\bf x)/\eps}){\bf g}({\bf x};\Lambda\bb{\sh(\bf x)/\eps})=\ord{\eps/\sh}$. If $\bf g\neq0$ then there may exist many $\lambda$ for which $\f({\bf x};\lambda)=0$, each generating a potential well, and hence creating many viable sliding modes near $\sh=0$. For large enough noise, the results of \cite{j13error,js2013} imply that the system eventually settles into the well corresponding to linear dependence on $\lambda$ (i.e. with ${\bf g}\equiv0$), leading to 
 \begin{equation}
 \begin{array}{llll}
\frac{d\langle{\bf x}\rangle}{dt}&\!\!\!\!=\!\!\!\!&\hf\bb{1+\lambda}\f ^+({\bf x})+\hf\bb{1-\lambda}\f ^-({\bf x})+(\lambda^2-1){\bf g}({\bf x};\lambda)+\ord{\kappa^2}&\!\mbox{for }\kappa<r(\eps)\;,\\
\frac{d\langle{\bf x}\rangle}{dt}&\!\!\!\!=\!\!\!\!&\hf\bb{1+\lambda}\f ^+({\bf x})+\hf\bb{1-\lambda}\f ^-({\bf x})+\ord{\kappa^2}&\!\mbox{for }\kappa>r(\eps)\;,
\end{array}
\end{equation}
for a function $r(\eps)$ whose form depends on ${\bf g}$, e.g. $r(\eps)=\sqrt{-\eps/\log\eps}$ for a friction example in \cite{js2013}. 

The counterintuitive outcome is that errors like noise can cause a system to behave more like a crude model (with linear switching) than a more refined one (with nonlinear switching), and hence discontinuous models owe their unreasonable effectiveness to unmodelled errors that wash out hidden effects of switching. But this washing out of nonlinearities is not universal. By analysing the ambiguity in how we treat the discontinuity we can quantify the effect of unmodelled errors, and estimate when they can be neglected.

\section{Forms \& origins of switching}\label{sec:bi}

%

In the literature on piecewise-smooth dynamical theory, much discussion is made of where discontinuous models are used, but little consideration is made of how discontinuities arise (though the question certainly occurs, e.g. in \cite{seidman96}). This is in part because the physical processes they model are often complicated or little understood, arising typically in engineering or biological or environmental contexts, and moreover they involve singular limits (as we shall see below), making the idea that a model lies `close to' some true system nontrivial. 
Let us therefore ask how discontinuities arise in the asymptotics of transition by means of various toy models, showing that the discontinuities that arise from ordinary and partial differential equations, from integral equations, or from heuristic sigmoid models, can be cast in a common form, namely \eref{Fasy}. 


Take a system that behaves like \eref{xasy}, 
where $\eps$ is a small positive constant that we ultimately set to $\eps=0$ to obtain a sharp transition. Let us assume the switch occurs due to a sudden transition in some extra variable $y$, scaled so that 
$y\sim \sign(\sh)$ for $|\sh|>\eps$, and propose that a complete model of the system can be written as
\begin{equation}\label{ds}
{\bf\dot x}=\F({\bf x};y)\qquad\mbox{such that}\qquad\f^\pm({\bf x})\equiv\F({\bf x},\pm1)\;.
\end{equation}
Our first task here is to show that broad classes all lead to asymptotic expressions of the form
\begin{equation}
y=\sign(\sh)+\ord{\eps/\sh}\quad\xrightarrow{\eps\rightarrow0}\quad\left\{\begin{array}{lll}+1&\rm if&\sh>0\;,\\-1&\rm if&\sh<0\;.\end{array}\right.
\end{equation}

\subsection{Ad hoc sigmoids}\label{sec:adhoc}

Piecewise-smooth dynamical theory has arisen chiefly to deal with situations where the precise laws of switching  are known. We should therefore begin our study by looking at the common empirical switching models, often ad hoc or based on incomplete physical intuition. 

One particular sigmoid function introduced by Hill \cite{hill} has become prevalent in biological models, and that is
$Z(z)=\frac{z^r}{z^r+\theta^r}$ 
for $z,\theta>0$, $r\in\mathbb N$. 
The function $Z(z)$ often represents the switching on/off of ligand binding or gene production in a larger model ${\bf\dot x}=\f({\bf x};y)$ of biological regulation. If we let $z=\theta e^\sh$ and $r=1/\eps$, for large argument the Hill function has an expansion 
$$y(\sh)=2Z(\theta e^\sh)-1=\sign(\sh)\cc{1-2e^{-|\sh|/\eps}+e^{-2|\sh|/\eps}+\ord{e^{-3|\sh|/\eps}}}\;.$$
In computation, commonly used sigmoids are the inverse or hyperbolic tangents, with expansions
$$\begin{array}{lcrcl}y(\sh)&=&\sfrac2\pi\arctan(\sh/\eps)&=&\sign(\sh)-\sfrac2\pi\cc{(\eps/\sh)+\ord{(\eps/\sh)^{3}}}\;,\\y(\sh)&=&\tanh(\sh/\eps)&=&\sign(\sh)\cc{1-2e^{-2|\sh|/\eps}+\ord{e^{-4|\sh|/\eps}}}\;,\end{array}$$ 
and one may expand various other sigmoids, like ${\sh/(\eps}{\sqrt{1+(\sh/\eps)^2})}$, 
in a similar way, with polynomially or exponentially small tails (i.e. $\ord{\eps/\sh}$ or $\ord{e^{-|\sh|/\eps}}$). 

Differentiable but non-analytic sigmoid functions are often used in theoretical approaches to smoothing discontinuities. An example is 
$$y(\sh)=\left\{\begin{array}{lll}r(-\sh)^{r(\sh)}-r(\sh)^{r(-\sh)}&\rm if&|\sh|<\eps\\\sign{\sh}&\rm if&|\sh|\ge\eps\end{array}\right.\qquad\qquad\qquad\qquad$$
where $r(\sh)=e^{2\eps/(\sh-\eps)}$. Its asymptotic form is rather messier than the examples above, but it is better behaved since its convergence to $\sign(\sh)$ is even faster, being given for $|\sh|<\eps$ by
$$y(\sh)=\sign(\sh)\cc{1-\frac{2e^{2/(|\sh/\eps|-1)}}{|\sh/\eps|+1}+\ord{e^{4/(|\sh/\eps|-1)}}-e^{2e^{-1}\cc{1+\ord{|\sh/\eps|-1}}/(|\sh/\eps|-1)}}\;.$$

In all cases the leading order term is made discontinuous by the presence of a $\sign(\sh)$, and the tails are small in $\sh/\eps$, being of order either $\ord{\eps/|\sh|}$, $\ord{e^{-|\sh|/\eps}}$, or $\ord{e^{1/(|\sh/\eps|-1)}}$. 
%
%
%

Friction (to be precise {dry}-friction) is a rich source of sigmoid switching models, with seemingly no limit to the different physical motivations and resulting laws. Yet the majority of arguments result in a dressed up $\sign$ function, a friction force $F(\sh)=\mu(\sh)\sign(\sh)$ where $\sh$ is the speed of motion along a rough surface and $\mu$ some smooth function, some including accelerative effects $F(\sh)=\mu(\sh,\dot \sh)\sign(\sh)$ or other nonlinearities to account for  `Stribeck' velocity or memory effects
(see e.g. \cite{wojewoda08,krim}); in almost all cases the $\sign$ function remains. 

\subsection{An ODE: Large-scale bistability, small-scale decay}\label{sec:log}

Let $y$ represent a population, for example, and consider a regulatory action that fixes $y$ to one constant value, $+1$, or another, $-1$, (up to some non-dimensionalization). During the transition the population might relax to a natural behaviour, decaying at a constant rate as $\dot y\sim-y$. 

Transitioning between steady states $y\sim \pm1$ for $|\sh|\gg\eps$ and relaxing as $\dot y\sim-y$ for $|\sh|\ll\eps$, for small $\eps$, is consistent with $(1-y^2)\sh=\eps(y+\dot y)$, and results in the system
\begin{equation}\label{log}
\begin{array}{rcl}
{\bf\dot x}&=&\f({\bf x};y)\;,\\
\eps\dot y&=&(1-y^2)\sh({\bf x})-\eps y\;.\\
\end{array}
\end{equation}
The quantity $\eps$ is small (the two $\eps$'s that appear here need not be the same, but for simplicity let us assume they are). 
Treating the $y$ system in \eref{log} as infinitely fast (for $\eps\rightarrow0$ so $\bf x$ is pseudo-static), its solution is easily found to be
\begin{equation}
y(t,\sh)=-(\eps/2\sh)+\alpha\tanh\bb{\alpha t\sh/\eps+k_0}\;, 
\end{equation}
where $\alpha=\sqrt{1+(\eps/2\sh)^2}$ and $k_0=\arctanh\bb{\frac{(\eps/2\sh)+y(0,\sh)}\alpha}$. This evolves on the fast timescale $t/\eps$ towards an attracting stationary state (where $\dot y=0>\partial\dot y/\partial y$), given by
$$\begin{array}{lll}
y_*(\sh)=-\eps/2\sh+\sign(\sh)\sqrt{1+(\eps/2\sh)^2}&\rm where&\frac{\partial\dot y}{\partial y}=-\sqrt{1+(2\sh/\eps)^2}\;.
\end{array}$$
For large $\sh$ the attractor sits close to either $+1$ or $-1$ depending on the sign of $\sh$. As $\sh$ passes through zero, $y_*(\sh)$ jumps rapidly (but continuously), and a series expansion for large $\sh/\eps$ reveals 
\begin{equation}
y_*(\sh)=\sign(\sh)-\frac{\eps}{2\sh}\cc{1-\frac{\eps}{4|\sh|}+\ord{(\eps/\sh)^3}}\;.
\end{equation} 
%
The asymptotic terms in the tail mitigate the transition in $|\sh|<\eps$, and everywhere else the variable $y$ relaxes to $y_*$ on a timescale $t=\ord\eps$, so we approximate $y\approx y_*=\sign(\sh)+\ord{\eps/\sh}$. 

\subsection{A PDE: Large-scale bistability, small scale dissipation}\label{sec:heat}

If instead $y$ represents a physical property like temperature, it might have both spatial and temporal variation that become significant during transition. 

For $|\sh|\ll\eps$ assume that $y$ satisfies the heat equation $y_{\sh\sh}\sim\eps\dot y$ for some small positive $\eps$, where $y_\sh$ denotes $\partial y/\partial \sh$. 
For $|\sh|\gg\eps$ assume asymptotes $y\sim\pm1$, implying $y_\sh\sim0$. 
This character is satisfied for example by the system 
$\frac \sh\eps y_\sh+y_{\sh\sh}-\eps\dot y=0$, giving overall
\begin{equation}\label{heat}
\begin{array}{rcl}
{\bf\dot x}&=&\f({\bf x};y)\;,\\
\eps^2\dot y&=&\sh({\bf x})y_\sh+\eps y_{\sh\sh}\;,
\end{array}
\end{equation}
The $y$ system evolves on a fast timescale $t/\eps^2$ to the slow subsystem $\sh({\bf x})y_\sh+\eps y_{\sh\sh}=0$, which has solutions $y=y_*(\sh)$ given by
\begin{equation}
y_*(\sh)=y_*(0)+y_{\sh*}(0)\sqrt{\pi\eps/2}\erf\sq{ \sh/{\sqrt{2\eps}}}\;,
\end{equation}
where $\erf$ denotes the standard error function \cite{as}. The asymptotes $y\rightarrow\pm1$ for large $\sh$ imply $y_*(0)=0$ and $y_{\sh*}(0)=\sqrt{2/\pi\eps}$. 
Solutions of the full system can be found in the form 
$y(t,\sh)=y_*(\sh)+e^{-t/\eps}Y(\sh)$. Substituting this into the partial differential equation for $y$ in \eref{heat} yields
\begin{eqnarray*}
0&=&\cc{\sh({\bf x})y_{\sh*}+\eps y_{\sh\sh*}}+e^{-t/\eps}\cc{\eps Y+\sh({\bf x})Y_\sh+\eps Y_{\sh\sh}}\;,
\end{eqnarray*}
(again treating $\bf x$ as pseudo-static for small $\eps$). 
The first bracket vanishes by the definition of $y_*$, the second gives an ordinary differential equation for $Y$ with solution
$$Y(\sh)=e^{-\sh^2/2\eps}\cc{y(0)\;_1F_1\sq{\sfrac{1-\eps}2,\frac12;\sfrac{\sh^2}{2\eps}}+\eps^{-1/2} y_\sh(0)\sh\;_1F_1\sq{1-\sfrac\eps2,\sfrac32;\sfrac{\sh^2}{2\eps}}}\;,$$
where $_1F_1$ is the Kummer confluent hypergeometric function \cite{as}. The exact functions are less interesting to us than their large variable asymptotics, given by
\begin{equation}
y_*(\sh)\sim \sign(\sh)-\frac{\sqrt{2\eps/\pi}}{\sh}e^{-\sh^2/2\eps}(1-\sqrt\eps/\sh+\ord{\eps/\sh^2})\;,
\end{equation}
(and for completeness,  
$Y(\sh)\sim \sqrt\pi\bb{\frac{y(0)}{\Gamma\sq{\frac{1-\eps}2}}+\frac{\sign(\sh)y_\sh(0)}{\sqrt2\Gamma\sq{1-\frac\eps2}}}\abs{\frac{\sqrt{2\eps}}v}^{\eps}+\ord{\eps/\sh^2}$). 

The function $Y(\sh)$ deviates from the sigmoid of $y_*(\sh)$ by an amount greatest near $\sh/\eps\approx0$ and decreasing inversely with $(\sh/\eps)^\eps$. Moreover this deviation disappears on the fast timescale $t/\eps$, so we approximate $y\approx y_*=\sign(\sh)+\ord{\sqrt\eps/\sh}$, similarly to \sref{sec:log} to leading order. 
%
%
(Evidently this system scales as $\sh/\sqrt\eps$ rather than $\sh/\eps$, a triviality fixed by replacing $\eps$ with $\eps^2$ in \eref{heat}).

\subsection{Integral turning points}\label{sec:int}

Lastly, let us turn from differential equations for $y$, to integrals. 
What follows is a very cursory description of a profound analytical phenomenon, for which we refer the reader to the literature as cited. 

\begin{figure}[h!]\centering\includegraphics[width=0.7\textwidth]{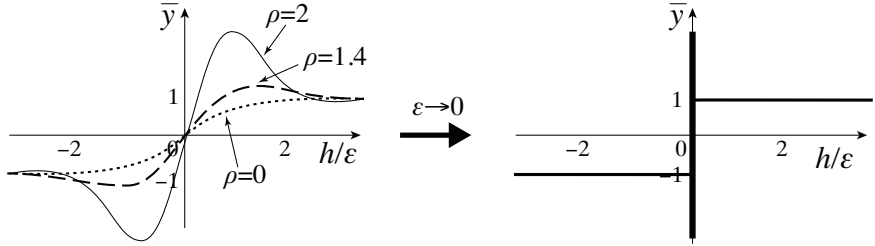}
\vspace{-0.3cm}\caption{\small\sf The graphs of $\bar y(\sh)$ for different values of $\rho$, which all limit to a $\sign$ function as $\eps\rightarrow0$. For $\rho>0$ the graph has peaks (multiple peaks for larger $\rho$), whose height is $\eps$-independent and therefore do not disappear as we shrink $\eps$, but merely get squashed into the region $|\sh|=\ord\eps$. }\label{fig:erfs}\end{figure}

First, as an example, take an integral over a simple Gaussian envelope $e^{-\hf\k^2}$, 
with a steady oscillation of wavelength $2\pi/\rho$, and an integration limit $h/\eps$, 
\begin{eqnarray}\label{erfdef}
y(\sh)&=&\sqrt{\frac2\pi}\int_{-\infty}^{\sh/\eps}\!d\k\; e^{-\hf\k^2}\cos(\rho\k)
\;.
\end{eqnarray}
Expanding this for large $h/\eps$ we obtain 
\begin{eqnarray}\label{erfc}
y(\sh)&\sim&e^{-\hf\rho^2}\bb{1+\sign(\sh)}-\frac{e^{-\sh^2/2\eps^2}}{\sqrt{\pi/2}}\cos{(\rho \sh/\eps)}\bb{\frac\eps{\sh}+\ord{(\eps/\sh)^3}}\;.
\end{eqnarray}
We can obviously now redefine $\bar y=e^{\hf\rho^2}y-1$ so that $\bar y=\sign h+...$ as in previous sections. 
Here $y$ is a simple sigmoid for $\rho=0$, but otherwise has peaks of height $|\bar y|\approx1+\sqrt{\frac2\pi}\frac{4\rho^3}{\pi^2}e^{-\pi^2/8\rho^2}$ at $h\approx\pm\eps\pi/2\rho$, illustrated in \fref{fig:erfs}. 
As we take the limit $\eps\rightarrow0$ for $\sh\neq0$, however, all graphs limit to $\bar y(\sh)=\sign(\sh)$ regardless of $\rho$, any peaks becoming squashed into the singular point $\sh=0$. 
%

The $\sign$ function here is the particularly well understood phenomenon of a {\it Stokes discontinuity} \cite{stokes64}. Their general role as a cause of discontinuities, associated with the rise and fall of large and small exponentials, requires innovative but not advanced application of complex variables, so a reasonable summary is warranted. 

More generally than \eref{erfdef}, say that $y$ is an integral of oscillations under an exponentially varying envelope, such as 
\begin{equation}\label{ints}
y(\alpha)=\int^{\alpha}_{-\infty}d\k\; a(\k)\;e^{\psi(\k)}\;.
\end{equation}
The term $a(\k)$ is taken to be slow (polynomially) varying, while the term $e^{\psi(\k)}$ is fast (exponentially) varying. 
This is typical when solving differential equations using Fourier or Laplace transforms, where typically $\psi(\k)=iu\k+\theta(\k)$ or $\psi(\k)=u\k+\theta(\k)$ respectively, where $u$ is a variable and $\k$ is its dual under the transform. (The fast varying term might not always be obvious, for example the transform of a high order polynomial $e^{u\k}[p(\k)]^N$ for large $N$ could be treated as an exponential $e^{\psi}$ with $\psi(\k)=u\k+N\log p(\k)$).
They can be analysed using stationary phase and steepest descent methods 
\cite{erdelyi,heading,dingle,bo99}. Care is needed in using them, but the principles are rather simple.


Assume that the integrand has a maximum at some point $\k$ along the integration path. If the integrand is oscillatory (when $\psi$ has an imaginary part), it will have many such maxima along the real line. But if the integrand is analytic then complex function theory allows us to deform the integration contour $(-\infty,\alpha]$ to anything of our choice in the complex plane of $\k$, provided it connects the point $-\infty$ to $\alpha$, and that we do not pass through infinities (e.g. poles of $ae^\psi$) in the process. If we could find a path $\op P$ along which the function was non-oscillating, and monotonic except perhaps for a maximum at some $\k_s$ where $\psi'(\k_s)=0$ (where $\psi'(\k)$ is the derivative with respect to complex $\k$), the approximation near this point would be
\begin{eqnarray}\label{stat1}
y&=\!&\int_{-\infty}^{\alpha}\!\! d\k\; a(\k)\;e^{\psi(\k)}\nonumber\\
&\approx\!&\int_{\op P}\!\! d\k\; \cc{a(\k_s)+(\k-\k_s)a'(\k_s)+...}\;e^{\cc{\;\psi(\k_s)\;+\;\hf(\k-\k_s)^2\psi''(\k_s)\;+\;...\;}}\nonumber\\
&\approx&a(\k_s)e^{\psi(\k_s)}\int_{\op P}\!\! d\k\;e^{\hf \;(\k-\k_s)^2\psi''(\k_s)}
\nonumber\\
&\approx&a(\k_s)e^{\psi(\k_s)}\frac{1}{\sqrt{-\psi''(\k_s)}}\int_{-\infty}^\infty\!\! du\;e^{-\hf u^2}
\;=\; a(\k_s)e^{\psi(\k_s)}\sqrt{\frac{2\pi}{-\psi''(\k_s)}}
\end{eqnarray}
to leading order. This is an incredibly simple, but also accurate, result, when properly used. In line 2 we just assume such an expansion is valid along a path $\op P$ (we will come back to this), and line 3 is just the leading order term. The clever bit is the simple substitution $u=(\k-\k_s)\sqrt{-\psi''(\k_s)}$ to obtain line 4, and this actually defines $\op P$ by demanding that $\op P$ transforms back to the real line. 

Some basic complex geometry makes all this work. Complex function theory tells us that the path $\op P$ we seek can indeed be found. 
By virtue of the Cauchy-Riemann equations, a path along which $\Im\psi=constant$ is also a steepest descent path of $\Re\psi$, so along such a path the function is non-oscillating (because the phase $\Im\psi$ is constant), and its magnitude $|e^{\psi}|=e^{\Re\psi}$ is exponentially fast varying (where $|e^{\psi}|$ is therefore exponentially). This only breaks down if the path encounters a maximum or minimum $\k_s$, where $\psi'(\k_s)=0$. That is exactly the point $\k_s$ which \eref{stat1} approximates about, integrating along the steepest descent path $\op P$, and the approximation is `exponentially good' because the integrand decays exponentially away from $\k_s$. 

We have neglected the endpoint $\alpha$. Because the integral is exponentially fast varying, the cutoff at the endpoint creates another exponentially strong maximum (or minimum, in which case we discard it), where typically $\psi'(\alpha)$ is non-vanishing. Approximating to leading order about $k=\alpha$ as above gives
\begin{eqnarray}\label{stat0}
y&\!=\!\!&\int_{-\infty}^{\alpha}\!\! d\k\; a(\k)\;e^{\psi(\k)}
\;\approx\; a(\alpha)e^{\psi(\alpha)}\int_{-\infty}^{\alpha}\!\!\!\!\!  d\k\;e^{ (\k-\alpha)\psi'(\alpha)}
\;\approx\;\frac{a(\alpha)e^{\psi(\alpha)}}{\psi'(\alpha)}\;.
\end{eqnarray}

So the endpoint, $\k=\alpha$, contributes to the integral if \eref{stat0} converges. The contribution of stationary point $\k_s$ is conditional, since it may or may not lie on the contour $\op P$, so we have 
\begin{eqnarray}\label{stat2}
y\approx-\frac{a(\alpha)e^{\psi(\alpha)}}{\psi'(\alpha)}+a(\k_s)e^{\psi(\k_s)}\sqrt{\frac{2\pi}{-\psi''(\k_s)}}\frac{1+\sign\sh}2\;.
\end{eqnarray}
The factor $\bb{1+\sign\sh}/2$ is a switch that turns on the stationary point contribution for $h>0$ if $\k_s\in\op P$, and turns it off for $h<0$ if $\k_s\notin\op P$. 
The transition between cases is a bifurcation in $\op P$ when the path connects $\k=\alpha$ directly to $\k=\k_s$, i.e. when $\Im\psi(0)=\Im\psi(\k_s)$ (since the path is a stationary phase contour). Typically we find up to a sign that
\begin{equation}
\sh=\Im\sq{\psi(0)-\psi(\k_s)}\;.
\end{equation}
In general there may be many stationary phase points $\k_{s1}$, $\k_{s2}$, ..., turned on and off at switching surfaces ({\it Stokes lines}) of the form $$\sh_{ i}=\Im\sq{\psi(\alpha)-\psi(\k_{si})}\qquad\mbox{or}\qquad\sh_{ij}=\Im\sq{\psi(\k_{si})-\psi(\k_{sj})}\;.$$ Finding the correct expansion \eref{stat2} requires inspection of the phase contours in the complex $\k$ plane, to find a path through the stationary points $\k_{si}$ and the endpoints of $\k\in(-\infty,\alpha]$, with $\op P$ permitted to pass through the `point at infinity' $|\k|\rightarrow\infty$, such that the integral converges. One may also calculate the higher order corrections, and a wealth of theory exists to assist, a good starting point is \cite{dingle}.

In the example \eref{erfdef}, the stationary point $\k_s=i\rho$ gives a contribution $2e^{-\hf\rho^2}$, and the endpoint $\k=\sh/\eps$ gives a contribution proportional to $e^{-\hf(\sh/\eps-i\rho)^2}$, and they have equal phase (they are both real) when $\sh=0$, providing the switching threshold. 

The point of all this is simply that, just as with our previous examples in this section, the discontinuity (the $\sign$ term) has again appeared as an inescapable part of the leading order behaviour \eref{stat2}, and remains there as we add higher orders in the tail of the series. The reader must pick apart the details to gain a fuller picture, but we have laid out the basics to illustrate how the sign function arises. 



\subsection{Return to the vector field}\label{sec:end}

Switching typically occurs when functions have different asymptotic behaviours on different domains, and this is what unites all of the examples above. The $\sign$ function affects the switch between different functional forms of $y$ that break down at $h=0$. 

The quantity $y$, which has a steady behaviour for almost all $h$, undergoes a sudden jump taking the form $y=\sign(h)+\alpha(\sh/\eps)\sum_{n=1}^\infty\beta_n(\eps/\sh)^n=\sign(h)+\ord{\eps/\sh}$. We then wish to model its effect on the dynamical system ${\bf\dot x}=\f({\bf x};y)$. In general $\f$ may have nonlinear dependence on $y$, as polynomials or trigonometric functions of $y$ for example, as in \eref{whirlnon} or \eref{vdp1}. The most we can then infer is that $\f$ takes a form $\f({\bf x};y(\sh))=\F({\bf x},h)$ as given by \eref{Fasy}. The consequences of that form are what we have presented already in this paper. 


\section{In closing}\label{sec:conc}

In \sref{sec:bi} we explored how discontinuities arise, not as crude modeling caricatures, but in the leading order of asymptotic expansions. Just as local expansions of differentiable functions yield linear terms, so asymptotic expansions of abrupt transitions yield discontinuities (characterized by the $\sign$ function here). They describe a switch in some unknown variable $y$, whose effect we then seek to understand on the bulk system ${\bf\dot x}=\f ({\bf x};y)$, using the methods of \srefs{sec:asydis}{sec:pws}. In practice, the origins of discontinuity explored in \sref{sec:bi} are often unknown, but we found them all to take a universal form, and we have shown how to express it in a manner that retains the asymptotic tails --- the ghosts of switching --- in the limit of a piecewise-smooth model. 


We have barely begun discovering the consequences of nonlinear switching for piecewise-smooth systems. The interaction of multiple switches, for example, opens up a vast world of attractors and bifurcations to be discovered. We have tried only to revisit the foundations of piecewise-smooth dynamics in a way that enables future study to embrace the ambiguity of the discontinuity, not to present a theory ready accomplished, and so many avenues are left to be explored in more rigour.

Discontinuities seem to be a symptom of interaction between incongruent objects or media, and the nature of such interactions is often difficult to model precisely. Whereas in some areas of physics we have a governing law, a wave or heat equation perhaps, to guide the transition or permit asymptotic matching, in many of the engineering and life science where discontinuous models are becoming increasingly prominent, we rely on much less perfect information. 

Piecewise-smooth dynamical theory attempts to address this, 
but we have seen that behaviour can be modelled that lies outside Filippov's simplest and most commonly adopted sliding theory. So how should we put nonlinear switching to use? Nonlinear terms offer more freedom to our switching paradigm, opportunities to model richer forms of dynamics that we have only begun to explore. The nonlinear terms may be matched to observations, or derived from physical principles if any are available, for example from a model like \eref{log}, \eref{heat}, or \eref{ints}, and in genetic regulation \cite{machina2013siads} or in friction such efforts are in progress. 

In addition we understand something of how sensitive a piecewise-smooth system is to its idealization of the switching as a discontinuity at a simple threshold. We can introduce a parameter $\eps$ characterizing stiffness if a switch is continuous (as in \sref{sec:adhoc}), and an amplitude $\kappa$ (or several $\kappa$'s) of discrete effects like noise, hysteresis, or time delay, again derived from physical principles or observations if possible, and in simulations the discretization step provides another $\kappa$ (as in \fref{fig:whirls}(b.ii)). The stiffness $\eps$ and unmodelled errors $\kappa$ compete, and in a $\kappa$ dominated system the nonlinear phenomena of hidden dynamics may be washed out, while they may flourish in a better behaved or better modelled (i.e. small $\kappa$) system.

That discontinuities yield strange dynamics is unsurprising, and the idea of `ghosts' left behind by approximation schemes is not new \cite{seidman96}. Perhaps more surprising is the extent to which we can characterise their effects in the piecewise-smooth framework. 
So what remains to be done? 
To the geometrical arsenal of singularities and bifurcations that we use to understand dynamical systems, we can add discontinuity-induced bifurcations \cite{bc08} and hidden attractors \cite{j15hbifs}. The task to classify these has a long way to go. 
Though it is not always made clear, many of the theoretical results in \cite{f88} (and hence to many works deriving from it) apply solely to the linear (or convex) combination found by assuming ${\bf g}\equiv0$. 
The nonlinear approach with ${\bf g}\neq0$ permits us to explore the different dynamics possible at the discontinuity, and thus to explore the many other systems that make up Filippov's full theory of differential inclusions. 
When nonsmooth systems do surprising things, we usually find we can make sense of them by extending our intuition for smooth systems to the switching layer, where, as in smooth systems, nonlinearity cannot be ignored.

Finally, there are currently no standard numerical simulation codes that can handle discontinuous systems with complete reliability, event detection being insufficient to take full account of all their singularities and issues of non-uniqueness (see e.g. \cite{j15hbifs,j15exit}). It is hoped that by capturing the ghosts of switching --- in the form of nonlinear discontinuity --- such codes may soon be developed.

%

\footnotesize

\bibliography{../../grazcat}
\bibliographystyle{plain} 

\end{document}